\documentclass[12pt]{amsart}
\usepackage{amsmath,amsfonts,amscd,amssymb}
\newtheorem{theorem}{Theorem}

\newtheorem{corollary}{Corollary}
\theoremstyle{remark}\newtheorem{remark}{Remark}
\theoremstyle{remark}

\begin{document}
\title[Enumeration of regular subgraphs]
{Tensor networks and the enumeration of regular subgraphs}

\author{Peter Zograf}
\address{Steklov Mathematical Institute, St.Petersburg, 191023 Russia}
\email{zograf@pdmi.ras.ru}

\thanks{Partially supported by  a RFBR grant.} 

\begin{abstract}
We propose a universal approach to a range of enumeration 
problems in  graphs. The key point is in contracting suitably chosen
symmetric tensors placed at the vertices of a graph along the edges. 
In particular, this leads to an algorithm that counts the number of 
$d$-regular subgraphs of an arbitrary graph including
the number of $d$-factors (previously we considered the case
$d=2$ with a special emphasis on the enumeration of Hamiltonian
cycles; cf. math.CO/0403339). We briefly discuss the problem 
of the computational complexity of this algorithm. \end{abstract}

\maketitle

\section{Introduction}

A {\em tensor network} is a collection of tensors placed at the vertices of a 
connected finite graph. The valency of each tensor must coincide with the 
degree of the corresponding vertex, and tensor indices are labelled by 
the half-edges incident to this vertex. Two tensors placed at adjacent
vertices can be contracted over the pair of indices labelling the two half-edges
of the connecting edge. This yields a tensor network with one tensor
less than the original one (the underlying graph obtained by contracting 
the corresponding edge has also one vertex and one edge less). When
contracting a loop, the tensor at its single vertex gets contracted over a
pair of indices, and the total number of tensors (or vertices) remains unchanged.
The total contraction along all edges of the graph gives a number (the graph
reduces then to a single vertex with no edges).

The idea of placing tensors at vertices of graphs and contracting them
along the edges traces back to Sylvester \cite{S}\footnote{We owe this
reference to A.~A.~Kirillov, Jr.}, who used this construction
for obtaining polynomial invariants of symmetric tensors. 
Much later Penrose \cite{P} reversed this idea
to computing graph invariants (like the number of 3-edge colorings of a planar 
3-regular graph)  in his novel (though not quite successful) approach to the 
4-color problem.  Later a number of other remarkable applications was found,
first of all  in the theory of Vassiliev knot invariants(see, e.~g.,~\cite{CDK}). Recently 
tensor network contraction was used to simulate quantum computation 
(cf. \cite{MS} and references therein).

Here we deal with contracting networks of specially chosen symmetric tensors.
The result of contraction gives a generating function for the numbers of regular
subgraphs of any type. In its simplest version this algorithm produces the
number of factors of any given degree. We also comment on the computational
complexity of this agorithm.

{\bf Acknowledgements.} We are grateful to S.~Cook and M.~Braverman for 
useful remarks.

\section{Networks of symmetric tensors and graph functions}\label{function}
Let $G$ be a finite graph (possibly with loops and multiple edges). 
The set of vertices 
of $G$ we denote by $V(G)=\{v_1, \dots ,v_n\}$,  where $n=|V(G)|$ 
is the total number of vertices, and
the set of edges of $G$ we denote by $E(G)$. For each vertex $v_i$
we denote by $d_i$ its degree (or valency), $i=1, \dots ,n$.
Then the number of edges of $G$ is given by 
$$|E(G)|=\frac{1}{2}\sum_{i=1}^n d_i.$$

Now let $\mathbb{F}$ be a field, and let $V\cong\mathbb{F}^r$ be a vector 
space of dimension $r$ over $\mathbb{F}$. Fix a symmetric bilinear form 
$B: V\otimes V\longrightarrow \mathbb{F}$. The graph $G$ 
together with the bilinear form $B$ define a multilinear form
\begin{equation}
B_G: V^{\otimes d_1}\otimes \dots \otimes V^{\otimes d_n}
\longrightarrow\mathbb{F},
\end{equation}
which is constructed as follows. At each vertex $v_i$ of $G$ we place  
$d_i$\nobreakdash-th tensor power $V^{\otimes d_i}$ of the vector space 
$V$, where the factors are labeled by the half-edges of $G$ incident to $v_i$. 
Each edge of $G$ defines a contraction of two copies of $V$ (corresponding 
to its two half-edges) by means of the bilinear
form $B$.  We obtain the multilinear form $B_G$ by performing 
such contractions over the set $E(G)$ of all edges of $G$.
Rigorously speaking, the multilinear form $B_G$ depends on the order of
half-edges at each vertex $v_i$, or, equvalently, on the order of factors in the
tensor power $V^{\otimes d_i}$. However, its restriction to $S^{d_1} V\otimes
\dots\otimes S^{d_n} V$, where $S^d V$ denotes the $d$\nobreakdash-th 
symmetric power of $V$, is defined uniquely. 

Now fix a sequence $\mathcal{A}=\{A_1,A_2, \dots \}$ of symmetric contravariant 
$d$\nobreakdash-valent tensors $A_d\in S^d V \subset V^{\otimes d}$. 
Here we consider tensor networks given by the triple $\{G,\,\mathcal{A},\,B\}$. We treat  
the tensor product $A_{d_1}\otimes \dots \otimes A_{d_n}$ as an element of 
$V^{\otimes d_1}\otimes \dots \otimes V^{\otimes d_n}$ 
and consider the element
\begin{equation}
\mathcal{F}_{\mathcal{A},\, B}(G)
=B_G (A_{d_1}\otimes \dots \otimes A_{d_i}) \in \mathbb{F}.
\end{equation}
In other words, $\mathcal{F}_{\mathcal{A},\, B}(G)$ is the result of contracting
the tensor network $\{G,\,\mathcal{A},\,B\}$ -- it is obtained by placing
a copy of $A_d$ at each vertex of $G$ of degree $d$ and contracting 
$\otimes_{i=1}^n A_{d_i}$
using $B$ over $|E(G)|$ pairs of indices corresponding to the edges of $G$.
Thus, to each pair $\mathcal{A},\, B$, where $\mathcal{A}$ is a sequence of
symmetric $d$-tensors ($d=1,2, \dots $) and $B$ is a
symmetric bilinear form, we associate an $\mathbb{F}$-valued mapping 
$\mathcal{F}_{\mathcal{A},\, B}$ on the set of isomorphism 
classes of graphs, or an  
$\mathbb{F}$-valued {\em graph function} in the terminology of \cite{T}.

\section{Enumeration of regular subgraphs}

Given a graph $G$, let $H$ be a $k$-regular subgraph of $G$
(not necessarily connected).
The {\em type} of the subgraph $H$ in $G$ is the partition 
$\lambda_H=[|V(H_1)|,\dots,|V(H_l)|]$
of the number $|V(H)|$ of vertices of $H$, 
where $H_1,\dots,H_l$ are the connected componets of
$H,\;l=l(H)$. The weight of partition $\lambda_H$ is $|\lambda_H|=|V(H)|$, and
the length is $l(\lambda_H)=l(H)$. For each partition $\lambda$ we define
a graph function $N_{k,\lambda}$ by 
$$N_{k,\lambda}(G)=\#\{H\subset G|H \text{ is $k$-regular and }\lambda_H=\lambda\},$$
i.e., $N_{k,\lambda}(G)$ is the number of $k$-regular subgraphs
of type $\lambda$ in $G$. Clearly, $N_{k,\lambda}(G)=0$ if $k>\max\{d_1,\dots,\,d_n\}$ or $|\lambda|>n$.

For $k$ a positive integer, denote by $p_k(x_1,x_2,\dots)=x_1^k+x_2^k+\dots$ 
the $k$\nobreakdash-th power sum in variables $x_1,x_2,\dots$.
Given a partition $\lambda=[k_1,\dots,k_l]$, we define a homogeneous
symmetric function $p_\lambda$ of degree $|\lambda|=k_1+\dots+k_l$ 
by the formula
$$p_\lambda(x_1,x_2,\dots)=\prod_{i=1}^l\;p_{k_i}(x_1,x_2,\dots).$$

We want to show that under a special choice of $\mathcal{A}$ and $B$ 
the graph function $\mathcal{F}_{\mathcal{A},\, B}$ defined in Section 
\ref{function} counts the number of $k$-regular subgraphs of any given type in graphs. 
We take $\mathbb{F}=\mathbb{C}$ and consider the standard 
coordinates in $V=\mathbb{C}^r$. In these coordinates 
the bilinear form $B$ is given by the identity $r\times r$ matrix $I_r$.
We define the tensors $A_{d,k}$ componentwise. For $d<k$ we put 
$$A_{d,k}^{i_1 \dots \,i_d}=\begin{cases} 0 & \text{if $(i_1 \dots \,i_d)\neq 
(r \dots r)$},\\ 
t & \text{if $(i_1 \dots \,i_d)= 
(r \dots r)$},\end{cases}$$
and for $d\geq k$ we put
$$A_{d,k}^{i_1 \dots \,i_d}=\begin{cases} x_i & \text{if $(i_1 \dots \,i_d)$ is a 
permutation of $(\underbrace{i\dots\,i}_k \,\underbrace{r \dots r}_{d-k})$},\\ 
& \hspace{2in} i=1, \dots ,r-1,\\ 
t & \text{if $(i_1 \dots \,i_d)=(r\,\dots\, r)$},\\
0 & \text{otherwise},\end{cases}$$
where $x_1,\dots, x_{r-1}$ and $t$ are arbitrary complex numbers.
The main result of this section is the following

\begin{theorem}
For $\mathcal{A}_k=\{A_{1,k},A_{2,k},\dots\}$ as above, 
the value of the graph function
$\mathcal{F}_{\mathcal{A}_k, \, I_r}$ on any graph $G$ is given 
by the formula
\begin{equation}
\mathcal{F}_{\mathcal{A}_k, \, I_r}(G)=\sum_{|\lambda|\leq n}
t^{n-|\lambda|}p_\lambda(x_1,\dots,x_{r-1})\, N_{k,\lambda}(G),
\label{poly}\end{equation}
where the sum is taken over the set of all partitions $\lambda$
of weight $|\lambda|\leq n=|V(G)|$.
\label{main}\end{theorem}

\begin{proof} 
We interprete the indices $1,\dots ,r$ as
colors of the half-edges of $G$. A product of $n$ components  
$A_{d_1,k}^{i_1 \dots\, i_{d_1}} \dots A_{d_n,k}^{i_{m-d_n+1} \dots\, i_m}$ 
(where $n=|V(G)|$ and $m=\sum_{j=1}^n d_j=2|E(G)|$) 
contributes non-trivially to $\mathcal{F}_{\mathcal{A}_k,\, I_r}(G)$ 
if and only if the colors agree on each edge of $G$ or, equivalently, 
if and only if for each edge the both indices that label two 
of its half-edges are the same. Thus, in this case 
the non-zero contributions are in one-to-one
correspondence with edge colorings of $G$ in $r$ colors 
with the following properties:

(i) an edge incident to a vertex of degree $<k$ has color $r$, and

(ii) at each vertex $v_j\in V(G)$ of degree $d_j\geq k$ 
exactly $k$ edges incident to it have  some  
color $i_j\in\{1, \dots ,r\}$, and the rest $d_i-k$ edges have color 
$r$ (if an edge makes a loop we count it twice).

The closure of the union of edges with colors $1, \dots, r-1$
is a $k$-regular subgraph $H$ in $G$,  and every connected component $H_j$ of $H$
is colored in one of the colors $i_j\in\{1, \dots ,r-1\}$. 
The contribution to $\mathcal{F}_{\mathcal{A}_k,\, I_r}(G)$ from this coloring
is $t^{n-|\lambda_H|}\prod_{j=1}^{l(H)} x_{i_j}^{|V(H_j)|}$, where 
$\lambda_H=[|V(H_1)|,\dots, |V(H_l)|]$ is the partition associated with $H$
and $l=l(H)$ is the number of connected components of $H$. 
Therefore, the contribution from all possible colorings of the subgraph 
$H$ is equal to 
$$t^{n-|\lambda_H|}\prod_{j=1}^{l(H)}\left(\sum_{i=1}^{r-1} x_i^{|V(H_j)|}\right)
=t^{n-|\lambda_H|}p_{\lambda_H}(x_1,\dots, x_{r-1}),$$
and summig up the contributions from all $k$-regular subgraphs in $G$ we get the 
assertion of the theorem.
\end{proof}

\begin{corollary}
The graph function $\mathcal{F}_{\mathcal{A}_k, \, I_r}$,
depending on $x_1,\dots, x_{r-1}$ and $t$ as parameters,  determines
the numbers $N_{k,\lambda}(G)$ uniquely for any graph
$G$ with $n\leq r-1$ vertices.
\end{corollary}

\begin{proof}
By Theorem \ref{main},  the graph function 
$\mathcal{F}_{\mathcal{A}_k, \, I_r}$  with values in $\mathbb{C}$
factors through the ring
$\mathbb{C}[x_1,\dots, x_{r-1}]^{S_{r-1}}$ of 
symmetric polynomials  
in $r-1$ independent variables $x_1,\dots, x_{r-1}$.
It is well known that the polynomials $p_k(x_1,\dots, x_{r-1}),\;
k=1,\dots, n$, are algebraically independent in
$\mathbb{C}[x_1,\dots, x_{r-1}]^{S_{r-1}}$ 
provided $n\leq r-1$. Therefore, in this case the graph function
$\mathcal{F}_{\mathcal{A}_k, \, I_r}(G)$ determines the coefficients
$N_{k,\lambda}(G)$ in (\ref{poly}) uniquely.
\end{proof}

\begin{remark}

Since the coefficients $N_{k,\lambda}(G)$ in (\ref{poly}) are
{\em non-negative} integers, we can uniquely find them out when
$$r\geq 1+\underset{H\subset G}{\max}\; l(H),$$
where $l(H)$ is the number of  connected componenets of $H$
and the maximum is taken over all $k$-regular subgraphs $H$ in $G$,
but we will not dwell on this issue here.
\label{remark}\end{remark}

Below are two special cases of Theorem \ref{main} of independent
interest.

\begin{corollary}
Put $r=2$, $x_1=1$ and $t=0$. Then for any graph $G$
the value $\mathcal{F}_{\mathcal{A}_k, \, I_2}(G)$ is
the number of {\em $k$-factors}, or spanning $k$-regular
subgraphs in $G$.
\label{H}\end{corollary}
\begin{proof}
By Theorem \ref{main}, 
$$\mathcal{F}_{\mathcal{A}_k, \, I_2}(G)=
x_1^n\,\sum_{|\lambda|=n}N_{k,\lambda}(G).$$
\end{proof}

In the simplest case when $k=1$, Corollary 1 counts the number 
of 1-factors, or {\em perfect mathcings} in $G$. 

The next statement concerns 
{\em connected} $k$-factors of $G$:

\begin{corollary}
Put $r=n+1$, $x_j=e^{2\pi\sqrt{-1}j/n},\; (j=1, \dots ,n)$ and $t=0$.
Then for any graph $G$ 
the number of connected $k$-factors
in $G$ is equal to
$\frac{1}{n}\,\mathcal{F}_{\mathcal{A}_k, \, I_{n+1}}(G)$, where
$n$ is the number of vertices of $G$.
\end{corollary}

\begin{proof}
In this case
$$p_k(x_1,\dots,x_n)=\begin{cases}
0, & k=1, \dots ,n-1,\\ n, & k=n.\end{cases}$$
By Theorem \ref{main} we have
$$\mathcal{F}_{\mathcal{A}_k, \, I_{n+1}}(G)=n\,,$$
where $[n]$ denotes the 1-element partition of $n=|V(G)|$
(by definition, $N_{k,[n]}(G)$ is the number of connected
$k$-factors).
\end{proof}

\section{Computational complexity}
The value
$\mathcal{F}_{\mathcal{A}_k, \, I_r}(G)$ can be effectively computed for any
graph $G$ as explained in Section \ref{function}. In general, this is a hard
computational problem. To make things less complicated let us take $r=2$. 
Then Corollary \ref{H} provides a relatively simple algorithm 
for counting the number of $k$-factors in graphs. (Note that even the case
$k=1$, or
enumeration of perfect matchings, is a well-known $\#P$-complete
problem.)

The computational complexity of the above algorithm depends on two main points:
\begin{enumerate}
\item the succession of tensor contractions along the edges of $G$ (or simply edge 
contractions); 
\item the utilization of the special form of tensors $A_{d,k}$; 
\end{enumerate}

What concerns the first point, the objective is to choose a sequence of edge contractions that keeps 
the maximal tensor valency (or vertex degree) as small as possible. This problem was solved in
\cite{MS}. Take a sequence of edge contractions reducing $G$ to a point. The {\em complexity}
of  this sequence is the maximum vertex degree during the contraction process. The {\em
contraction complexity} $cc(G)$ is the minimal complexity over all sequences of edge contractions.
A nice result of \cite{MS} states that $cc(G)=tw(G^*)+1$, where $G^*$ is the linear graph
of $G$, and $tw(G^*)$ is the {\em treewidth} of $G^*$.  Though computing the treewidth of a general
graph is $NP$-hard, a tree decomposition of $G^*$ of width $O(tw(G^*))$ can be obtained in time
$|V(G^*)|^{O(1)}\exp(O(tw(G^*))$ \cite{RS}. Moreover, given a tree decomposition of $G^*$ of
width $w$, a sequence of edge contractions in $G$ of complexity not greater than $w+1$ can
be found in polynomial time. Thus, the results of \cite{MS} give a rough estimate of the 
computational complexity of our algorithm by $|V(G)|^{O(1)}\exp(O(tw(G^*)))$. So, for graphs
whose line graphs have bounded treewidth our algorithm works quite fast. Moreover, known
constructions of almost optimal tree-decompositions (cf. \cite{RS}) give almost optimal edge
contraction sequences.

As for the second point mentioned above, the situation is less clear at the moment.
Note that the estimate of \cite{MS} is rather general and applies to arbitrary tensor networks.
In our case it can be considerably sharper, since the tensor networks we use are of very
special types.  In particular,  Sylvester's theory of decomposing symmetric tensors into 
sums of symmetric powers of  vectors may appear quite useful  in improving the efficiency 
of our algorithm.\footnote{This observation belongs to V.~N.~Vassiliev.} However, this is still a 
work in progress.


\begin{thebibliography}{66}

\bibitem{CDK} S.~V.~Chmutov, S.~V.~Duzhin and A.~I.~Kaishev, 
{\it The algebra of 3-graphs},
Trans. Steklov Math. Inst. {\bf 221}, 1998, 157-186. \par
\bibitem{MS} I.~Markov, Y.~Shi, {\it Simulating quantum computation by
contracting tensor networks}, Preprint arXiv: quant-ph/0511069.\par
\bibitem{P} R.~Penrose, {\it Applications of negative dimensional tensors},
Combinatorial mathematics and its applications (ed.  D.~J.~A.~Welsh),
Academic Press, 1971, 221-244. \par
\bibitem{RS} N.~Robertson, P.~D.~Seymour, {\it Graph minors. X. Obstructions to
tree-decompositions}, J.~Combin. Theory Ser.~B, {\bf 52}, 1991, 153-190.\par
\bibitem{S} J.~Sylvester, {\it Sur les covariants irr\'{e}ductibles du quantic binare de 
huiti\`{e}me order. Collected mathematical papers of J.~J.~Sylvester, vol. III}, 
Cambridge Univ. Press, 1909, 481-488.\par
\bibitem{T} W.~T.~Tutte, {\it Graph theory}, Addison-Wesley, 1984.\par
\bibitem{Z} P.~Zograf, {\it The enumeration of edge colorings and Hamiltonian
cycles by means of symmetric tensors}, Preprint arXiv: math.CO/0403339. \par

\end{thebibliography}
\end{document}